\documentclass[12pt]{article}
\usepackage{mathrsfs}
\usepackage{amscd}
\usepackage{amsmath,amsfonts,amssymb,amscd}
\usepackage{indentfirst,graphics,epsfig,psfrag}
\input{epsf}
\usepackage{ifpdf}
\usepackage{enumerate}
\usepackage{appendix}
\usepackage{enumerate}

\usepackage{lineno}

\makeatletter \@addtoreset{figure}{section} \makeatother
\makeatletter
\long\def\@makecaption#1#2{%
   \vskip 10\p@
   \setbox\@tempboxa\hbox{{#1}\ \ #2}%
   \ifdim \wd\@tempboxa >\hsize
       {#1}\ \ #2\par
   \else
       \hbox to\hsize{\hfil\box\@tempboxa\hfil}%
   \fi}
\makeatother

\newtheorem{thm}{Theorem}[section]
\newtheorem{cor}[thm]{Corollary}
\newtheorem{lem}[thm]{Lemma}

\newtheorem{obs}[thm]{Observation}
\newtheorem{pro}[thm]{Proposition}

\newcommand{\qed}{{\hfill\rule{3pt}{7pt}}}

\def\qed{\hfill \rule{4pt}{7pt}}

\begin{document}
\title{Rainbow connection numbers of complementary graphs\footnote{Supported by NSFC.}}
\author{
\small  Xueliang Li, Yuefang Sun\\
\small Center for Combinatorics and LPMC-TJKLC\\
\small Nankai University, Tianjin 300071, P.R. China\\
\small E-mails: lxl@nankai.edu.cn, syf@cfc.nankai.edu.cn
 }
\date{}
\maketitle
\begin{abstract}

A path in an edge-colored graph, where adjacent edges may be colored
the same, is a rainbow path if no two edges of it are colored the
same. A nontrivial connected graph $G$ is rainbow connected if there
is a rainbow path connecting any two vertices, and the rainbow
connection number of $G$, denoted by $rc(G)$, is the minimum number
of colors that are needed in order to make $G$ rainbow connected. In
this paper, we provide a new approach to investigate the rainbow
connection number of a graph $G$ according to some constraints to
its complement graph $\overline{G}$. We first derive that for a
connected graph $G$, if $\overline{G}$ does not belong to the
following two cases: $(i)$~$diam(\overline{G})=2,3$,
$(ii)~\overline{G}$ contains exactly two connected components and
one of them is trivial, then $rc(G)\leq 4$, where $diam(G)$ is the
diameter of $G$. Examples are given to show that this bound is best
possible. Next we derive that for a connected graph $G$,
if $\overline{G}$ is triangle-free, then $rc(G)\leq 6$. \\[2mm]
\end{abstract}

\section{Introduction}

All graphs in this paper are finite, undirected and simple. Let $G$
be a nontrivial connected graph on which an edge-coloring $c:
E(G)\rightarrow \{1,2,\cdots,n\}$, $n\in \mathbb{N}$, is defined,
where adjacent edges may be colored the same. A path is $rainbow$ if
no two edges of it are colored the same. An edge-coloring graph $G$
is $rainbow ~connected$ if any two vertices are connected by a
rainbow path. Clearly, if a graph is rainbow connected, it must be
connected. Conversely, any connected graph has a trivial
edge-coloring that makes it rainbow connected; just color each edge
with a distinct color. Thus, we define the
$rainbow~connection~number$ of a connected graph $G$, denoted by
$rc(G)$, as the smallest number of colors that are needed in order
to make $G$ rainbow connected. If $G'$ is a connected spanning
subgraph of $G$, then $rc(G)\leq rc(G')$. Chartrand et al. obtained
that $rc(G)=1$ if and only if $G$ is complete, and that $rc(G)=m$ if
and only if $G$ is a tree, as well as that a cycle with $k> 3$
vertices has rainbow connection number $\lceil \frac{k}{2}\rceil$, a
triangle has rainbow connection number 1 (\cite{Chartrand 1}). Also
notice that, clearly, $rc(G)\geq diam(G)$ where $diam(G)$ denotes
the diameter of $G$. In an edge-colored graph $G$, we use $c(e)$ to
denote the color of an edge $e$, and for a subgraph $H$ of $G$,
$c(H)$ denotes the set of colors of edges in $H$. We use~$V(G)$,
$E(G)$ for the set of vertices and edges of $G$, respectively. For
any subset $X$ of $V(G)$, denote $G[X]$ as the subgraph induced by
$X$, and $E[X]$ the edge set of $G[X]$; For a set $S$, $|S|$ denotes
the cardinality of $S$. As usual, $P_n$ is a path on $n$ vertices.
For a connected graph $G$, the $distance$ between two vertices $u$
and $v$ in $G$, denoted by $dist(u,v)$, is the length of a shortest
path between them in $G$. The $eccentricity$ of a vertex $v$ in $G$
is defined as $ecc_G(v)=\max_{x\in V(G)}{dist(v,x)}$. We follow the
notation and terminology of \cite{Bondy}.

In this paper, we provide a new approach to investigate the rainbow
connection number of a graph $G$ according to some constraints to
its complement graph $\overline{G}$. We give two sufficient
conditions to guarantee that $rc(G)$ is bounded by a constant.

One of our main results is:
\begin{thm}\label{thm5} For a connected
graph $G$, if $\overline{G}$ does not belong to the following two
cases: $(i)$~$diam(\overline{G})=2,3$, $(ii)~\overline{G}$ contains
exactly two connected components and one of them is trivial, then
$rc(G)\leq 4$. Furthermore, this bound is best possible.\qed
\end{thm}

For the remaining cases, $rc(G)$ can be very large as discussed in
Section 4. So we add a constraint, i.e., we let $\overline{G}$ be
triangle-free. Then $G$ is claw-free, and we can derive our next
main result:
\begin{thm}\label{thm9} For a connected graph $G$, if $\overline{G}$
is triangle-free, then $rc(G)\leq 6$.\qed
\end{thm}

\section{Preliminaries}

We now give a necessary condition for an edge-colored graph to be
rainbow connected. If $G$ is rainbow connected under some
edge-coloring, then for any two cut edges (if exist) $e_1=u_1u_2$,
$e_1=v_1v_2$, there must exist some $1\leq i,j\leq 2$, such that any
$u_i-v_j$ path must contain edge $e_1,e_2$. So we have:
\begin{obs}\label{ob1}
If $G$ is rainbow connected under some edge-coloring, $e_1$ and
$e_2$ are any two cut edges, then $$c(e_1)\neq c(e_2).$$\qed
\end{obs}

The following lemma will be useful in our discussion.
\begin{lem}[\cite{Y. Caro}]\label{thm101}
If $G$ is a connected graph and $H_1,\cdots,H_k$ is a partition of
the vertex set of $G$ into connected subgraphs, then $rc(G)\leq
k-1+\sum_{i=1}^k{rc(H_i)}$.\qed
\end{lem}

In \cite{Chartrand 1}, the authors derived the precise values of the
rainbow connection numbers of complete bipartite graph
$K_{s,t}(2\leq s\leq t)$ and complete $k$-partite graph ($k\geq 3$).

\begin{thm}[\cite{Chartrand 1}]\label{thm102} For integers $s$ and $t$ with $2\leq s\leq
t$, $$rc(K_{s,t})=\min\{{\lceil \sqrt[s]{t} \rceil}, 4\}.$$\qed
\end{thm}

\begin{thm}[\cite{Chartrand 1}]\label{thm103}
Let $G=K_{n_1,n_2,\ldots,n_k}$ be a complete $k$-partite graph,
where $k\geq 3$ and $n_1\leq n_2\leq \ldots \leq n_k$ such that
$s=\sum_{i=1}^{k-1}{n_i}$ and $t=n_k$. Then

\[
 rc(G)=\left\{
   \begin{array}{ll}
     1 &\mbox {if~$n_k=1$,}\\
     2 &\mbox {if~$n_k\geq 2$~and~$s>t$,}\\
     \min\{\lceil \sqrt[s]{t}\rceil,3\} &\mbox {if~$s\leq t$.}
   \end{array}
   \right.
\]\qed

\end{thm}

From the above two theorems, we know that $rc(K_{s,t})\leq 4$ for
any $s,t \geq 2$ and $rc(G)\leq 3$ where $G$ is a complete
$k$-partite graph with $k\geq 3$.

We now introduce a definition from \cite{Chandran}, A dominating set
$D$ in a graph $G$ is called a $two-way~ dominating~set$ if every
pendant vertex of $G$ is included in $D$. In addition, if $G[D]$ is
connected, we call $D$ a $connected~two-way~ dominating~set$. Note
that if $\delta(G)\geq 2$, then every $(connected)$ dominating set
in $G$ is a (connected) two-way dominating set. We also need the
following result.

\begin{thm}[\cite{Chandran}]\label{thm104} If $D$ is a connected
two-way dominating set in a graph $G$, then $rc(G)\leq
rc(G[D])+3$.\qed
\end{thm}

\section{Proof of Theorem \ref{thm5}}

We first investigate the rainbow connection numbers of connected
complement graphs of graphs with diameter at least 4.
\begin{thm}\label{thm1} Let $G$ be a connected graph with $diam(G)\geq 4$. If
$\overline{G}$ is connected, then $rc(\overline{G})\leq 4$.
\end{thm}
\begin{pf} We choose a vertex $x$ with $ecc_{G}(x)=diam(G)=d\geq 4$.
Let $N_G^i(x)=\{v: dist(x,v)=i \}$ where $0\leq i\leq d$. So
$N_G^0(x)=\{x\}, N_G^1(x)=N_G(x)$ as usual. Then $\bigcup_{0\leq
i\leq d}{N_G^i(x)}$ is a vertex partition of $V(G)$ with
$|N_G^i(x)|=n_i$. Let $A=\bigcup_{i~is~even}{N_G^i(x)}$,
$B=\bigcup_{i~is~odd}{N_G^i(x)}$. For example, see Figure
\ref{figure1}, a graph with $diam(G)=4$.
\begin{figure}[!hbpt]
\begin{center}
\includegraphics[scale=1.000000]{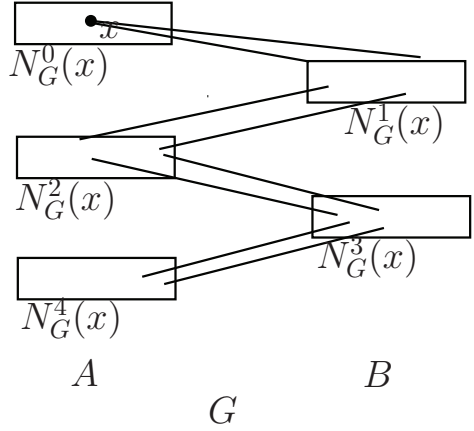}
\end{center}
\caption{Graph for the example with $d=4$.}\label{figure1}
\end{figure}

So, if $d=2k(k\geq 2)$ then $A=\bigcup_{0\leq i\leq
d~is~even}{N_G^i(x)}$, $B=\bigcup_{1\leq i\leq
d-1~is~odd}{N_G^i(x)}$; if $d=2k+1(k\geq 2)$ then $A=\bigcup_{0\leq
i\leq d-1~is~even}{N_G^i(x)}$, $B=\bigcup_{1\leq i\leq
d~is~odd}{N_G^i(x)}$. Then by the definition of complement graphs,
we know that $\overline{G}[A](\overline{G}[B])$ contains a spanning
complete $k_1$-partite subgraph (complete $k_2$-partite subgraph)
where $k_1=\lceil \frac{d+1}{2} \rceil (k_2=\lceil \frac{d}{2}
\rceil)$. For example, see Figure \ref{figure1}, $\overline{G}[A]$
contains a spanning complete tripartite subgraph $K_{n_0,n_2,n_4}$,
$\overline{G}[B]$ contains a spanning complete bipartite subgraph
$K_{n_1,n_3}$.

\textbf{Case 1.} $d\geq 5$. Then $k_1,k_2\geq 3$. From Theorem
\ref{thm103}, we have $rc(\overline{G}[A]), rc(\overline{G}[B])\leq
3$.

We now give $\overline{G}$ an edge-coloring as follows: we first
give the subgraph $\overline{G}[A]$ a rainbow edge-coloring using
three colors; then give the subgraph $\overline{G}[B]$ a rainbow
edge-coloring using the same colors as that of the subgraph
$\overline{G}[A]$; next we give a fresh color to all edges between
the subgraph $\overline{G}[A]$ and the subgraph $\overline{G}[B]$.

We will show that this coloring is rainbow. It suffices to show that
for any $u\in \overline{G}[A]$, $v\in \overline{G}[B]$, there is a
rainbow path connecting them in $\overline{G}$. We first choose an
edge $uv_1$ where $v_1\in \overline{G}[B]$ (it must exist, without
loss of generality, we assume $u\in N_G^2(x)$, then $u$ is adjacent
to all vertices in $N_G^5(x)$). Then by adding a rainbow $v_1-v$
path in $\overline{G}[B]$, we obtain our desired path. So
$rc(\overline{G})\leq 4$ in this case.

\textbf{Case 2.} $d= 4$, that is, $A={N_G^0(x)}\cup {N_G^2(x)}\cup
{N_G^4(x)}$, $B={N_G^1(x)}\cup {N_G^3(x)}$. So
$\overline{G}[A](\overline{G}[B])$ contains a spanning complete
3-partite subgraph $K_{n_0,n_2,n_4}$ (complete bipartite subgraph
$K_{n_1,n_3}$). So, from Theorem \ref{thm103} we have
$rc(\overline{G}[A])\leq 3$.

\textbf{Subcase 2.1.} $n_1,n_3\geq 2$. Since now $\overline{G}[B]$
contains a spanning complete bipartite subgraph $K_{n_1,n_3}$, from
Theorem \ref{thm102} we have $rc(\overline{G}[B])\leq 4$.

We now give $\overline{G}$ an edge-coloring as follows: we first
give the subgraph $\overline{G}[B]$ a rainbow edge-coloring using
four colors, say $a,b,c,d$; then give the subgraph $\overline{G}[A]$
a rainbow edge-coloring using colors $a,b,c$; next we give the color
$d$ to all edges between the subgraph $\overline{G}[A]$ and the
$\overline{G}[B]$.

We will show that this coloring is rainbow. It suffices to show that
for any $u\in \overline{G}[A]$, $v\in \overline{G}[B]$, there is a
rainbow path connecting them in $\overline{G}$. We first choose an
edge $vu_1$ where $u_1\in \overline{G}[A]$ (it must exist, without
loss of generality, we assume $v\in N_G^1(x)$, then $v$ is adjacent
to all vertices in $N_G^4(x)$). Then by adding a rainbow $u_1-u$
path in $\overline{G}[B]$, we obtain our desired path. So
$rc(\overline{G})\leq 4$ in this case.

\textbf{Subcase 2.2.} At least one of $n_1,n_3$ is 1, say $n_1=1$.

We now give $\overline{G}$ an edge-coloring as follows: we give the
edges between $N_G^0(x)$ and $N_G^4(x)$ a color $a$; give the edges
between $N_G^0(x)$ and $N_G^2(x)$ a new color $b$; give the edges
between $N_G^2(x)$ and $N_G^4(x)$ a new color $c$; give the edges
between $N_G^1(x)$ and $N_G^4(x)$ a new color $d$; give the edges
between $N_G^0(x)$ and $N_G^3(x)$ the color $b$; give the edges
between $N_G^1(x)$ and $N_G^3(x)$ the color $c$.

We will show that this coloring is rainbow. We only need to show
that there is a rainbow path connecting two vertices $u,v\in
N_G^3(x)$, the remaining cases are easy. Let $P:= u,x,x_1,x_2,v$
where $x_1\in N_G^4(x), x_2\in N_G^1(x)$. Clearly, it is rainbow. So
$rc(\overline{G})\leq 4$ in this case. \qed
\end{pf}

With a similar argument to that of Theorem \ref{thm1}, we have:
\begin{pro}\label{thm3}
If $G$ is a tree but not a star, then $rc(\overline{G})\leq 3$.
\end{pro}
\begin{pf} It is easy to show that if $G$ is a tree but not a
star, then $\overline{G}$ is connected. We now use the same
terminology as in the argument of Theorem \ref{thm1}. Note that $A$
and $B$ are independent sets in $G$ (consider the BFS-tree of $G$).
So, $\overline{G}[A]$ and $\overline{G}[B]$ are two disjoint cliques
in $\overline{G}$. Then by Lemma \ref{thm101} we have
$rc(\overline{G})\leq 3$. \qed
\end{pf}

Theorem \ref{thm1} is equivalent to the following result.
\begin{thm}\label{thm2} For a connected graph $G$, if $\overline{G}$
is connected and $diam(\overline{G})\geq 4$, then $rc(G)\leq 4$.\qed
\end{thm}

If $G$ is a graph with $h\geq 2$ connected components, then
$\overline{G}$ contains a complete $h$-partite spanning subgraph,
and so we have
\begin{pro}\label{thm4} If $G$ is a graph with $h\geq 2$ connected components
$G_i$ and $n_i'=n(G_i)(1\leq i\leq h)$, then $rc(\overline{G})\leq
rc(K_{n_1',\cdots,n_h'})$.\qed
\end{pro}

\textbf{Proof of Theorem \ref{thm5}.} If $\overline{G}$ is
connected, since $diam(\overline{G})\neq 2,3$ and clearly
$diam(\overline{G})\neq 1$, from Theorem \ref{thm2} we have
$rc(G)\leq 4$. If $\overline{G}$ is disconnected, since by the
assumption, it has either at least three connected components or
exactly two nontrivial components, then from Theorems \ref{thm102}
and \ref{thm103} and Proposition \ref{thm4} we have $rc(G)\leq 4$.

Let $\overline{G}$ contain two connected components, one is a clique
with $s\geq 2$ vertices, the other is a clique with $t\geq 3^s+1$
vertices. We have $G=K_{s,t}$, then from Theorem \ref{thm102},
$rc(G)=\min\{{\lceil \sqrt[s]{t} \rceil}, 4\}=4$, and so the bound
is best possible. \qed

\section{Proof of Theorem \ref{thm9}}

For the remaining cases, since the complement of $\overline{G}$ is
$G$ itself, we need to investigate $rc(\overline{G})$ in two cases:
$(i)$~$diam(G)=2,3$, $(ii)~G$ contains two connected components and
one of them is trivial. We first give some discussion about the case
$diam(G)=3$. We use the same terminology as that of Theorem
\ref{thm1}.

\begin{thm}\label{thm6} For a vertex $x$ of $G$ satisfying
$ecc_G(x)=diam(G)=3$, we have $rc(\overline{G})\leq 5$ for the three
cases $(i)~n_1=n_2=n_3=1$, $(ii)~n_1,n_2=1,n_3\geq 2$, and
$(iii)~n_2=1,n_1,n_3\geq 2$. For the remaining cases,
$rc(\overline{G})$ may be very large. Furthermore, if $G$ is
triangle-free and $\overline{G}$ is connected, then
$rc(\overline{G})\leq 5$.
\end{thm}
\begin{pf} If $n_1=n_2=n_3=1$, then $G$ is a 4-path $P_4$, and so
$rc(\overline{G})=3$. Thus, we could consider the following three
cases.

\textbf{Case 1.} Two of $n_1,n_2,n_3$ are equal to 1.

\textbf{Subcase 1.1.} $n_1,n_2=1$. Then it is easy to show that the
subgraph $\overline{G}[N_G^0(x)\cup N_G^1(x)\cup N_G^3(x)]$ contains
a bipartite spanning subgraph $K_{2,n_3}$, and so from Lemma
\ref{thm101} and Theorem \ref{thm102} we have $rc(\overline{G})\leq
rc(K_{2,n_3})+1\leq 5$.

\textbf{Subcase 1.2.} $n_1,n_3=1$. Let $n_2'=|\{v\in N_G^2(x):
deg_{\overline{G}}(v)=1 \}|$. Then there are $n_2'$ cut edges in
$\overline{G}$, and so from Observation \ref{ob1} we have
$rc(\overline{G})\geq n_2'$.

Furthermore, if $G$ is triangle-free, then $N_G^2(x)$ is a stable
set in $G$, and so a clique in $\overline{G}$, and thus from Lemma
\ref{thm101} we have $rc(\overline{G})\leq 4$.

\textbf{Subcase 1.3.} $n_2,n_3=1$. With a similar argument to that
of \textbf{Subcase 1.2}, we have $rc(\overline{G})\geq n_1'$ where
$n_1'=|\{v\in N_G^1(x): deg_{\overline{G}}(v)=1 \}|$.

Furthermore, if $G$ is triangle-free, then $N_G^1(x)$ is a stable
set in $G$, and so a clique in $\overline{G}$, and thus from Lemma
\ref{thm101} we have $rc(\overline{G})\leq 4$.

\textbf{Case 2.} One of $n_1,n_2,n_3$ is equal to 1.

\textbf{Subcase 2.1.} $n_1=1$. With a similar argument to that of
\textbf{Subcase 1.2}, we have $rc(\overline{G})\geq n_2'$ where
$n_2'=|\{v\in N_G^2(x): deg_{\overline{G}}(v)=1 \}|$.

Furthermore, if $G$ is triangle-free, then $N_G^2(x)$ is a stable
set in $G$, and so a clique in $\overline{G}$. In $\overline{G}$,
the subgraph $\overline{G}[N_G^0(x)\cup N_G^1(x)\cup N_G^3(x)]$
contains a spanning bipartite subgraph $K_{2,n_3}$. So from Theorem
\ref{thm102}, it needs at most four colors to rainbow it; we then
give a new color to the edges between $x$ and $N_G^2(x)$. Clearly,
this coloring is rainbow and we have $rc(\overline{G})\leq 5$.

\textbf{Subcase 2.2.} $n_2=1$. Then it is easy to show that the
subgraph $\overline{G}[N_G^0(x)\cup N_G^1(x)\cup N_G^3(x)]$ contains
a spanning bipartite subgraph $K_{1+n_1,n_3}$. So from Lemma
\ref{thm101} and Theorem \ref{thm102}, we have $rc(\overline{G})\leq
rc(K_{1+n_1,n_3})+1\leq 5$.

\textbf{Subcase 2.3.} $n_3=1$. Let $N_G^3(x)=\{u\}$. With a similar
argument to that of \textbf{Subcase 1.2}, we have
$rc(\overline{G})\geq n_1'+n_2'$ where $n_i'=|\{v\in N_G^i(x):
deg_{\overline{G}}(v)=1 \}|$ with $i=1,2$.

Furthermore, if $G$ is triangle-free, then $N_G^1(x)$ is a stable
set in $G$, and so a clique in $\overline{G}$. Let $V_u$ be the set
of vertices of $N_G^2(x)$ which are adjacent to $u$ in $G$. So $V_u$
is a stable set in $G$ and a clique in $\overline{G}$. We now give
$\overline{G}$ an edge-coloring: We give the edges of the complete
graph $\overline{G}[N_G^1(x)\cup\{u\}]$ a color $a$; give the edge
$xu$ a new color $b$, give the edges (they may not exist, but now
$N_G^2(x)=V_u$ is a clique and the procedure is easy) between $u$
and $N_G^2(x)\backslash V_u$ a new color $c$; the edges between $x$
and $N_G^2(x)$ a new color $d$. It is easy to check that the
coloring is rainbow and $rc(\overline{G})\leq 4$ in this case.

\textbf{Case 3.} $n_1,n_2,n_3\geq 2$. With a similar argument to
that of \textbf{Subcase 1.2}, we have $rc(\overline{G})\geq n_2'$
where $n_2'=|\{v\in N_G^2(x): deg_{\overline{G}}(v)=1 \}|$.

Furthermore, if $G$ is triangle-free, then $N_G^1(x)$ is a stable
set in $G$, and so a clique in $\overline{G}$. If every vertex in
$N_G^3(x)$ is adjacent to all vertices of $N_G^2(x)$ in $G$, then
both $N_G^2(x)$ and $N_G^3(x)$ are stable sets in $G$, and so
cliques in $\overline{G}$, since $G$ is triangle-free. Then in
$\overline{G}$, $\overline{G}[N_G^0(x)\cup N_G^3(x)]$,
$\overline{G}[N_G^2(x)]$, $\overline{G}[N_G^1(x)]$ are complete
graphs. So from Lemma \ref{thm101} we have $rc(\overline{G})\leq
rc(\overline{G}[N_G^0(x)\cup
N_G^3(x)])+rc(\overline{G}[N_G^2(x)])+rc(\overline{G}[N_G^1(x)])+2=5$.
Thus we choose a vertex $u\in N_G^3(x)$ with $V_u \neq \emptyset,
N_G^2(x)$, where $V_u$ denotes the set of neighbors of $u$ in
$N_G^2(x)$ in $G$, and so it is a stable set in $G$ and a clique in
$\overline{G}$.

We now give $\overline{G}$ an edge-coloring: We give a new color $a$
to the edges of $\overline{G}[N_G^1(x)]$; for every vertex $w$ of
$N_G^3(x)\backslash \{u\}$, since $w$ is adjacent to all vertices of
$N_G^1(x)$ in $\overline{G}$, we give a new color $b$ to an edge
between $w$ and $N_G^1(x)$, give a new color $c$ to the remaining
edges between $w$ and $N_G^1(x)$; give color $a$ to the edges
between $x$ and $N_G^3(x)\backslash \{u\}$; give the edge $xu$ a new
color $d$; give a new color $e$ to the edges between $x$ and
$N_G^2(x)$; give the color $b$ to the edges between $u$ and
$N_G^2(x)\backslash V_u$. It is easy to check that the above
coloring is rainbow and $rc(\overline{G})\leq 5$ in this case.

From the above discussion, we know that $rc(\overline{G})\leq 5$ for
the three cases $(i)~n_1=n_2=n_3=1$, $(ii)~n_1,n_2=1,n_3\geq 2$, and
$(iii)~n_2=1,n_1,n_3\geq 2$. For the remaining cases,
$rc(\overline{G})$ can be very large if $n_i'(i=1,2)$ is
sufficiently large. Furthermore, if $G$ is triangle-free, then
$rc(\overline{G})\leq 5$.\qed
\end{pf}

The following corollary clearly holds.

\begin{cor}\label{thm7} For a connected graph $G$, if $\overline{G}$
is triangle-free and $diam(\overline{G})=3$, then $rc(G)\leq 5$.\qed
\end{cor}

For a graph $G$ with $diam(G)=2$, let $x$ be a vertex satisfying
$ecc_G(x)=diam(G)$. Then, the two cases: $(i)~n_1=n_2=n_3=1$ and
$(ii)~n_1=1,n_2\geq 2$ do not hold, since in both cases
$\overline{G}$ are disconnected and $rc(\overline{G})~are~
undefined$. For the remaining two cases, that is, $n_1\geq 2,n_2=1$,
$n_1,n_2\geq 2$, with a similar argument to that of Theorem
\ref{thm6}, we have $rc(\overline{G})\geq n_1'$,
$rc(\overline{G})\geq n_2'$, respectively. So $rc(\overline{G})$ can
be very large if $n_i'(i=1,2)$ is sufficiently large. So we add an
additional constraint, i.e., we let $G$ be triangle-free.

\begin{pro}\label{thm10} Let $G$ be a triangle-free graph with
$diam(G)=2$. If $\overline{G}$ is connected, then
$rc(\overline{G})\leq 5$.
\end{pro}
\begin{pf} We choose a vertex $x$ with $ecc_G(x)=diam(G)=2$, and we use
the same terminology as that of Theorem \ref{thm1}. By the above
discussion, we only need to consider the following two cases.

\textbf{Case 1.} $n_1\geq 2,n_2=1$. Since $G$ is triangle-free,
$N_G^1(x)$ is a stable set in $G$ and so a clique in $\overline{G}$.
Thus, $rc(\overline{G})\leq 3$.

\textbf{Case 2.} $n_1,n_2\geq 2$. Since $G$ is triangle-free,
$N_G^1(x)$ is a stable set in $G$ and so a clique in $\overline{G}$.
Since $\overline{G}$ is connected, there exist $u\in N_G^1(x), v\in
N_G^2(x)$ such that $uv\in E(\overline{G})$.

If there exists some vertex $w\in N_G^2(x)$ with $deg_G(w)=n-2$,
then $w$ is adjacent to the remaining vertices except $x$ in $G$.
Since $diam(G)=2$, there exists $w_1w_2 \in E(G)$ with $w_1\in
N_G^1(x), w_2\neq w \in N_G^2(x)$. So $\{w,w_1,w_2\}$ is a triangle
in $G$, this produces a contradiction. So $deg_G(w)<n-2$ for all
$w\in N_G^2(x)$, and $deg_{\overline{G}}(w)\geq 2$ for all $w\in
N_G^2(x)$. Let $D=\{x,v\}\cup N_G^1(x)$. Then $D$ is a connected
two-way dominated set in $\overline{G}$. So from Theorem
\ref{thm104}, we have $rc(\overline{G})\leq
rc(\overline{G}[D])+3\leq 5$. \qed
\end{pf}

If $G$ contains two connected components, say $G_1, G_2$. Let
$n_1'=|\{v\in G_2: deg_G(v)=n-2 \}|$. Then in $\overline{G}$, there
are $n_1'$ pendant vertices and so there are $n_1'$ cut edges. From
Observation \ref{ob1}, we have $rc(\overline{G})\geq n_1'$. So in
this case, $rc(\overline{G})$ can be very large if $n_1'$ is
sufficiently large. So we also add an additional constraint, i.e.,
we let $G$ be triangle-free.

\begin{pro}\label{thm8} If $G$ is triangle-free and contains two connected
components one of which is trivial, then $rc(\overline{G})\leq 6$.
\end{pro}
\begin{pf} Suppose that $G$ contains two components, one is trivial,
the other is not trivial. Since $G$ is triangle-free, then
$\overline{G}$ is claw-free. Let $u$ be the isolated vertex in $G$,
so it is adjacent to any other vertex in $\overline{G}$, and so
$diam(\overline{G})=2$. We will consider two cases according to the
value of $\delta_{\overline{G}}$ where $\delta_{\overline{G}}$
denotes the minimum degree of $\overline{G}$.

\textbf{Case 1.} $\delta_{\overline{G}}=1$. Let
$deg_{\overline{G}}(v_1)=\delta_{\overline{G}}$ and $v_1v_2\in
\overline{G}$ ($v_2=u$). Since $\overline{G}$ is claw-free, the
subgraph $\overline{G}[V\backslash \{v_1\}]$ is a complete graph, so
$rc(\overline{G})=2$.

\textbf{Case 2.} $\delta_{\overline{G}}\geq 2$. Let
$deg_{\overline{G}}(v_1)=\delta_{\overline{G}}$. Then $u\in
N_{\overline{G}}^1(v_1)$ and is adjacent to any other vertex in
$\overline{G}$. So the subgraph $\overline{G}[D]$ contains a
spanning bipartite subgraph $K_{2,\delta_{\overline{G}}-1}$ where
$D=\{v_1\}\cup N_{\overline{G}}^1(v_1)$. Clearly, $D$ is a connected
two-way dominating set. We give the edge $uv_1$ a color $a$, give
the edges between $v_1$ and $N_{\overline{G}}^1(v_1)\backslash
\{u\}$ a new color $b$, give the edges between $u$ and
$N_{\overline{G}}^1(v_1)\backslash \{u\}$ a new color $c$. It is
easy to check that this coloring is rainbow. From Theorems
\ref{thm104} and \ref{thm102}, we have $rc(\overline{G})\leq
rc(\overline{G}[D])+3\leq 6$.  \qed
\end{pf}

From Theorem \ref{thm5}, Corollary \ref{thm7}, and Propositions
\ref{thm10} and \ref{thm8}, our next main result can be derived.

\textbf{Proof of Theorem \ref{thm9}.} We consider two cases:

\textbf{Case 1.} $\overline{G}$ is connected. The result holds for
the case $diam(\overline{G})\geq 4$ from Theorem \ref{thm5}, the
case $diam(\overline{G})= 3$ from Corollary \ref{thm7} and the case
$diam(\overline{G})= 2$ from Proposition \ref{thm10}.

\textbf{Case 2.} $\overline{G}$ is disconnected. The result holds
for the case that $\overline{G}$ contains two connected components
with one of them trivial from Proposition \ref{thm8}, and holds for
the remaining case from Theorem \ref{thm5}.\qed


\begin{thebibliography}{99}%参考文献的写法
\small \setlength{\itemsep}{-.8mm}
\bibitem{Bondy} J.A. Bondy, U.S.R. Murty,
{\it Graph Theory}, GTM 244, Springer, 2008.

\small \setlength{\itemsep}{-.8mm}
 \bibitem{Y. Caro} Y. Caro, A. Lev, Y. Roditty, Z. Tuza, R. Yuster,
 $On$ $rainbow$ $connection$, Electron. J. Combin. {\bf 15} (2008), R57.

\small \setlength{\itemsep}{-.8mm}
\bibitem{Chandran} L.S. Chandran, A. Das, D. Rajendraprasad,
N.M. Varma, $Rainbow$ $connection$ $number$ $and$ $connected$
$dominating$ $sets$, Arxiv preprint arXiv:1010.2296v1 [math.CO]
(2010).

\small \setlength{\itemsep}{-.8mm}
\bibitem{Chartrand 1} G. Chartrand, G.L. Johns, K.A. McKeon,
P. Zhang, $Rainbow$ $connection$ $in$ $graphs$, Math. Bohem. {\bf
133}(2008) 85-98.


\end{thebibliography}
\end{document}